\documentclass[12pt]{amsart}

\newcommand{\sm}{\setminus}

\newcommand{\inv}{^{-1}}

\newcommand{\wt}{\widetilde}

\newcommand{\PP}{{\mathbb P}}

\newcommand{\C}{{\mathbb C}}

\renewcommand{\phi}{\varphi}

\newcommand{\Ga}{\Gamma}

\newcommand{\ep}{\varepsilon}

\newenvironment{rem}{\medskip\noindent{\it Remark:\/} }{\medskip}

\title{New examples of hyperbolic octic surfaces
in $\PP^3$}

\author{Bernard Shiffman}
\address{Department of Mathematics, Johns Hopkins University, Baltimore,
MD 21218, USA}
\email{shiffman@math.jhu.edu}

\author{Mikhail Zaidenberg}
\address{Universit{\'e}
Grenoble I, Institut Fourier, UMR 5582 CNRS-UJF, BP 74,
38402 St.\ Martin
d'H{\`e}res c{\'e}dex, France}
\email{zaidenbe@ujf-grenoble.fr}

\date{June 24, 2003}

\begin{document}

\begin{abstract} We show that a general small deformation of the union
of two general cones in $\PP^3$ of degree $\ge 4$ is Kobayashi
hyperbolic.  Hence we obtain new examples of  hyperbolic
surfaces in
$\PP^3$ of any given degree $d\ge 8$.
\end{abstract}

\maketitle

It was shown by Clemens \cite{Cl} that  a  very general surface
$X_d$ of degree $d\ge 5$ in $\PP^3$ has no rational curves;  G.
Xu \cite{Xu} showed that $X_d$  also has no elliptic curves (and
in fact has no curves of genus $\le 2$), i.e.\ $X_d$ is
algebraically hyperbolic. According to the Kobayashi Conjecture,
$X_d$ must even be Kobayashi hyperbolic, and hence does not
possess non-constant entire curves $\C\to X_d$. The latter
property is known to be open in the Hausdorff topology on the
projective space of degree $d$ surfaces \cite{Za}, and it does
hold for a very general surface of degree at least $15$
\cite{DEG,El,MQ}.

Examples of hyperbolic surfaces in $\PP^3$ have been given by many
authors; see the references in our previous papers
\cite{ShZa1,ShZa2}, where more examples are given. So far,
the minimal degree of known examples is $8$; the first family of
 examples  of degree 8 hyperbolic surfaces in $\PP^3$  was found
by Fujimoto \cite{Fu} and independently by Duval \cite{Du}. In
\cite{ShZa2}, we introduced a deformation method, which we used to
construct a new degree 8 hyperbolic surface. In this note, we use a
simple form of our deformation method to construct another degree $8$
example, which is  a deformation of the union of two quartic cones.
Actually, our construction provides examples in any degree $d\ge 8$.

It follows from an observation by Mumford and Bogomolov, proved
in \cite{MM}, that every surface in $\PP^3$ of degree at most 4
contains rational or elliptic curves. However, it remains unknown
whether there exist hyperbolic surfaces in $\PP^3$ in the
remaining degrees $d=5,6,7$.

\medskip To describe our examples, we consider an algebraic curve $C$
in a plane $H\subset \PP^3$. We let  $\langle C,p\rangle$ denote
the cone formed by the union of  lines through a fixed point
$p\in\PP^3\sm H$ and points of $C$.  By a cone in $\PP^3$, we mean
a cone  of the form $X=\langle C,p\rangle$.  If $C'=X\cap H'$,
where $H'$ is an arbitrary plane not passing through $p$, then we
also have $X=\langle C',p\rangle$. We observe that $\deg X=\deg
C$.

\medskip\noindent {\bf Theorem.} {\it For $m,n\ge 4$, a general small
deformation of the union $X=X'\cup X''$ of two general cones in $\PP^3$
of degrees
$m$ and $n$,
respectively, is a hyperbolic surface of degree $m+n$. }

\begin{proof} Let $X=X'\cup X''\subset \PP^3$ be the union of two general
cones of respective degrees $m,n\ge 4$. We choose coordinates
$(z_1:z_2:z_3:z_4)\in\PP^3$ so that $X',X''$ are cones through the points
$a=(0:0:0:1)$ and $b=(0:0:1:0)$ respectively.  We consider the planes
$H'=\{z_4=0\},\ H''=\{z_3=0\}$ in $\PP^3$, and we write
\begin{eqnarray*}F_1=X'\cap H'=\{f_1(z_1,z_2,z_3)=0,\ z_4=0\}\;,\\
F_2=X''\cap H'' =\{f_2(z_1,z_2,z_4)=0,\ z_3=0\}\;,\end{eqnarray*}
where $f_1,f_2$ are  general homogeneous polynomials of degree
$m,n$ respectively. As we noted above, we have $X'=\langle F_1,
a\rangle$, $X''=\langle F_2, b\rangle$; hence, $X$ is the surface
of degree $m+n$ with equation
$$f_1(z_1,z_2,z_3)f_2(z_1,z_2,z_4)=0\;,\qquad (z_1:z_2:z_3:z_4)\in\PP^3.$$

We assume that $a\not\in X''$ and $b\not\in X'$, i.e.\ $f_1(0,0,1)\neq 0$
and $f_2(0,0,1)\neq 0$.
Let $$\pi_0:\PP^3\ \dashrightarrow \PP^1\;,\qquad
(z_1:z_2:z_3:z_4)\mapsto (z_1:z_2)$$ be the projection from the line
$z_1=z_2=0$. We further assume that $F_1$ and $F_2$ are smooth and
that each fiber of
$\pi_0|F_1$ and of $\pi_0|F_2$ has at least 3 distinct points. For
example, if
$m=n=4$, this will be the case whenever $(0:0:1)$ does not lie on any of
the bitangents or inflection tangent lines of $\{f_1=0\}$ or $\{f_2=0\}$.

We follow the deformation method of our paper \cite{ShZa2}. Let
$X_\infty=\{f_{\infty}=0\}$ be a general surface of degree $m+n$
in $\PP^3$, and let
$$X_t=\{f_1(z_1,z_2,z_3)f_2(z_1,z_2,z_4)+t\,f_{\infty}(z_1,z_2,z_3,z_4)=0\}\qquad
(t\in\C)\;.$$ We claim that $X_t$ is hyperbolic for sufficiently small
$t\neq 0$.  Suppose on the contrary that $X_{t_n}$ is not hyperbolic for a
sequence $t_n\to 0$. Then by Brody's Theorem \cite{Br}, there exists a
sequence
$\phi_n:\C\to X_{t_n}$ of entire holomorphic curves such that
$$\|\phi_n'(0)\| =\sup_{w\in \C} \|\phi_n'(w)\| =1,\qquad
n=1,2,\dots\,$$ where the norm is measured with respect to the
Fubini-Study metric in $\PP^3$.  Hence after passing to a
subsequence, we can assume that $\phi_n$ converges to a
nonconstant entire curve $\phi:\C\to X$.

Since $X=X'\cup X''$, we may suppose without loss of generality
that $\phi(\C)\subset X'$. Consider the projection from $a$,
$$\pi_a:\PP^3\dashrightarrow\PP^2\,,\qquad (z_1:z_2:z_3:z_4)\mapsto
(z_1:z_2:z_3).$$ Then $f_1\circ \pi_a\circ\phi=0$; i.e,
$\pi_a\circ\phi(\C)\subset \{f_1=0\}\approx F_1$.  Since $F_1$ is
hyperbolic (it has genus $\ge 3$), it follows that $\pi_a\circ\phi$ is
constant, and hence
$\phi(\C)$ is contained in a projective line of the form
$$\langle
p,a\rangle=\{(p_1s_0:p_2s_0:p_3s_0:s_1)\in\PP^3\ |\
(s_0:s_1)\in\PP^1\}\;,$$ where $p=(p_1:p_2:p_3:0)\in F_1$. We
notice that $(p_1,p_2)\neq (0,0)$ by the hypothesis that
$f_1(0,0,1)\neq 0$.

Let $\Gamma=X'\cap X''$ denote the double curve of $X$, and
suppose that $q\in\Gamma\cap\langle p,a\rangle=X''\cap\langle
p,a\rangle$. Recalling that $b\not\in X'$, we see that $q$ is of
the form $q=(p_1:p_2:p_3:s)$ where $f_2(p_1:p_2:s)=0$.  Thus we
have a bijection
$$\Gamma\cap\langle p,a\rangle
\buildrel {\approx}\over \to
\pi_0\inv(p_1:p_2)\cap F_2\;,\qquad (p_1:p_2:p_3:s) \mapsto (p_1:p_2:0:s)\;.$$ For
general $x=(p_1:p_2)$ the set $\pi_0\inv(x)\cap F_2$ contains
$n\ge 4$ distinct points, and by our assumption, it contains at
least 3 distinct points for all $x\in\PP^1$. Hence $\Gamma\cap\langle p,a\rangle$
contains at least 3 distinct points for all $p\in F_1$, and contains $n$ points for
general $p\in F_1$.

\medskip \noindent {\it Claim:\ \ $\phi(\C) \subset \langle p,a\rangle \sm
(\Gamma \sm X_\infty).$}

\medskip \noindent {\it Proof of the claim:\/}
 Suppose on the contrary that
$$\phi(w_0)=(\zeta_1:\zeta_2:\zeta_3:\zeta_4)\in
\Gamma
\sm X_\infty$$ for some $w_0\in\C$.
Let
$\Delta$ be a small disk about
$w_0$ such that
$\phi(\bar\Delta)\cap X_\infty=\emptyset$.  After shrinking $\Delta$ if
necessary, we can lift the maps $\phi_n|\bar\Delta$ via the
projection $\pi:\C^4\sm\{0\}\to\PP^3$ to maps
$\wt\phi_n:\bar\Delta\to
\C^4$ such that $$\wt\phi_n\to \wt\phi,\qquad
\pi\circ\wt\phi=\phi|\bar\Delta.$$
(Simply choose $j$ with $\zeta_j\neq 0$ and let $(\wt \phi_n)_j\equiv 1$.
Note that by our hypothesis that
$a,b\not\in\Gamma$, we can choose $j=1$ or 2.)

Let $n$ be  sufficiently large so that $\phi_n(\bar\Delta)$ does
not meet $X_\infty$. Then $f_{\infty}\circ \wt\phi_n$ does not
vanish on $\bar\Delta$. Since $\phi_n(\bar\Delta)\subset X_t$, it
then follows from the equation for $X_t$  that $f_2\circ
\wt\phi_n$ cannot vanish on $\bar\Delta$ (where we write
$f_2(z_1,z_2,z_3,z_4)=f_2(z_1,z_2,z_4)$). On the other hand, since
$\phi(w_0)\in X''$, we have $f_2\circ \wt\phi(w_0)=0$. It then
follows from Hurwitz's Theorem that $f_2\circ \wt\phi \equiv 0$,
i.e.\ $\phi(\Delta)\subset X''$. Then $\phi$ is constant since
$\phi(\Delta)$ lies in the finite set $X''\cap\langle
p,a\rangle$, a contradiction.  This verifies the claim.

\medskip
We now assume that, for all $p\in F_1$, the set $\langle
p,a\rangle \cap (\Gamma \sm X_\infty)$ contains at least 3 points,
or in other words, the finite set $X_{\infty}\cap \Gamma$ does not
contain 2 distinct points of $\langle p,a\rangle$, and does not
contain any of the points $\Gamma\cap\langle p,a\rangle$ for the
special values of $p$ where $\Gamma\cap\langle p,a\rangle$
consists of only 3 points. Similarly, we make the same assumption
for $F_2$. To show that this assumption holds for general
$X_\infty$, we consider the branched cover
$$\pi_\Gamma:=\pi_0|\Gamma:\Gamma\to \PP^1.$$ General fibers of $\pi_\Ga$
contain $mn$ distinct points. It suffices to show that a general
$X_\infty$
\begin{enumerate}
\item[i)] does not contain 2
distinct points of any fiber of $\pi_\Ga$ (i.e.\ $\pi_0|(\Ga\cap X_\infty)$
is injective), and
\item[ii)] does not contain any of the
points of the special fibers with fewer than $mn$ points.
\end{enumerate}
Since the totality of points in (ii) is finite, (ii) certainly
holds for general $X_\infty$.  It then suffices to show (i) for
the nonspecial fibers. Since  $\pi_\Ga$ is nonbranched at the
points of the nonspecial fibers, these points are smooth points
of $\Gamma$, and hence by Bertini's theorem, a general divisor
$X_\infty$ intersects $\Gamma$ transversally at these points. Now
suppose that $X_\infty=\{f_{\infty}=0\}$ intersects $\Ga$
transversally and does not intersect the special fibers, and
furthermore $\pi_0(\Ga\cap X_\infty)$ has maximal cardinality
among such $X_\infty$.  If (i) does not hold, then we can write
$\Ga\cap X_\infty=\{q^1,q^2,\dots,q^{(m+n)mn}\}$, where
$\pi_0(q^1)=\pi_0(q^2)$. Choose a divisor $Y=\{h=0\}$ of degree
$m+n$ containing the point $q^1$ but not $q^2$, and let
$X_\infty^\ep=\{f_{\infty}+\ep h=0\}$. For small $\ep$, we let
$q^j_\ep$ denote the point of $\Ga\cap X_\infty^\ep$ close to
$q^j$. (These points are well defined and the maps $\ep\mapsto
q^j_\ep$ are continuous for small $\ep$, since by the
transversality assumption, $f_{\infty}|\Ga$ has only simple
zeros.) Then $q^1_\ep=q^1$ and $q^2_\ep \neq q^2$ for small
$\ep\neq 0$. Hence for $\ep$ sufficiently small,
$\pi_0(q^2_\ep)\neq \pi_0(q^2)=\pi_0(q^1_\ep)$ and $\#[
\pi_0(\Ga\cap X_\infty^\ep)]
> \#[ \pi_0(\Ga\cap X_\infty)]$, a contradiction.

Thus, $\langle p,a\rangle \cap
(\Gamma \sm X_\infty)$ contains at least 3 points, for all $p\in F_1$
(for general $X_\infty$). Since $\phi(\C) \subset \langle p,a\rangle \sm
(\Gamma \sm X_\infty)$, $\phi$ is constant by Picard's
Theorem, which is a contradiction.
\end{proof}

\begin{rem} For an alternative construction of surfaces with hyperbolic
deformations, we let
$F=\{f=0\}$ and $G=\{g=0\}$ be two general plane curves of
degrees $m\ge 4$ and $n\ge 2$, respectively. We suppose that the
projective line $z_0=0$ meets $F$ ($G$, respectively)
transversally at $m$ ($n$, respectively) distinct points
$\{a_1,\ldots,a_m\}$ ($\{b_1,\ldots,b_n\}$, respectively). We
then consider the following cones in $\PP^4$ (with coordinates
$(z_0:\ldots :z_4)$) over these curves:
$$Y_1:=\langle F, u\rangle=\{f(z_0,z_1,z_2)=0\}\qquad \mbox{and}\qquad
Y_2:=\langle G,v\rangle=\{g(z_0,z_3,z_4)=0\}\;,$$
where the vertex sets are the skew projective lines
$$u:=\{z_0=z_1=z_2=0\}\qquad \mbox{and}\qquad v:=\{z_0=z_3=z_4=0\}\,.$$
We let $Y:=Y_1\cap Y_2\,.$ Thus $Y$ is an irreducible complete
intersection surface in $\PP^4$ of degree $mn$. It has $m+n$
singular points $\{A_1,\ldots,A_m\}=v\cap X$ of multiplicity $n$
and $\{B_1,\ldots,B_n\}=u\cap Y$ of multiplicity $m$ and no
further singularities. Indeed, the hyperplane section $Y\cap
H_{\infty}$, where $H_{\infty}:=\{z_0=0\}\simeq \PP^3$, is the
union of $mn$ distinct projective lines $l_{jk}:=\langle
A_jB_k\rangle\,\,(j=1,\ldots,m,\,k=1,\ldots,n)$, $n$ lines through each
point $A_j$ and $m$ through each point $B_k$.

Then $Y$ is birational to the direct product $F\times G$. Indeed,
it is obtained by blowing up the $mn$ points $c_{jk}:=a_j\times
b_k\in F\times G\,\,(j=1,\ldots,m,\,k=1,\ldots,n)$, and then
blowing down the proper transforms of $m$ vertical generators
$a_j\times G$ and $n$ horizontal generators $F\times b_k$ to the
singular points $A_j\in X$ and $B_k\in
X,\,\,j=1,\ldots,m,\,k=1,\ldots,n$, respectively.

We let now $\pi: \PP^4\dashrightarrow H_{\infty}\simeq\PP^3$ be a
general projection with center $P_0=(1:0:0:0)\notin Y\cup
H_{\infty}$, and we let $Z:=\pi(Y)\subset\PP^3$ (with the coordinates
$(z_1:z_2:z_3:z_4)$). Then $Z$  is given by
the resultant $r:=\mbox{res}_{z_0} (f(z_0,z_1,z_2),
g(z_0,z_3,z_4))$.

One can easily check in the same way as above that a general
small deformation of $Z$ is a hyperbolic surface in $\PP^3$ of
degree $mn\ge 8$.

The degenerate case $g=(z_0-z_3)(z_0-z_4)$  gives again the
union of two cones $X=X'\cup X''$ as in  the above theorem for the case
$f_1=f_2=f$.
\end{rem}

\end{document}